\documentclass[12pt]{amsart}
\usepackage{amscd}
\usepackage{amssymb}
\usepackage{a4wide}
\usepackage{amstext}
\usepackage{amsthm}

\makeatletter
\newcommand*{\rom}[1]{\expandafter\@slowromancap\romannumeral #1@}
\newcommand{\specialnumber}[1]{%
  \def\tagform@##1{\maketag@@@{(\ignorespaces##1\unskip\@@italiccorr#1)}}%
}
\newcommand{\specialeqref}[2]{\begingroup
  \def\tagform@##1{\maketag@@@{(\ignorespaces##1\unskip\@@italiccorr#2)}}%
  \eqref{#1}\endgroup}
\makeatother

\newtheorem{theorem}{Theorem}

\newtheorem{corollary}[theorem]{Corollary}

\newtheorem{lemma}[theorem]{Lemma}

\newtheorem{proposition}[theorem]{Proposition}

\newcommand{\R}{{\mathbb R}}\newcommand{\T}{{\mathbb T}}

\newcommand{\Z}{{\mathbb Z}}
\newcommand{\N}{{\mathbb N}} 
\newcommand{\La}{\Lambda}\newcommand{\la}{\lambda}\newcommand{\si}{\sigma}

\begin{document}

\title{On the Duality between Sampling and Interpolation}
\author{Alexander Olevskii  and  Alexander Ulanovskii}

\begin{abstract} We present a simple method based on the stability and duality  of the properties of sampling and interpolation, which allows one to substantially simplify the proofs of some classical results.

\end{abstract}

\maketitle


\section{Introduction}
\subsection{Sampling and Interpolation Problems}
A set $\La\subset\R$ is called uniformly discrete (u.d.) if \begin{equation}\label{sepc}d(\La):=\inf_{\la,\la'\in\La,\la\ne\la'}|\la-\la'|>0.\end{equation}
The constant $d(\La)$ is called the separation constant for $\La$.

In signal processing, sampling means the reduction of a continuous signal to a u.d. set $\La$. Two classical problems are of a great importance:

\medskip

(i) The {\it Sampling Problem} concerns with a reconstruction of continuous signals from their samples;

\medskip
(ii) The {\it Interpolation Problem} concerns with the possibility of representation of discrete functions on $\La$ as the traces of continuous signals.

\medskip
The classical papers of A. Beurling, J.-P.  Kahane and H. Landau deal with these problems in the Paley--Wiener and Bernstein settings.

The aim of this paper is to present a simple method which stresses the duality between sampling and interpolation problems, and allows to substantially simplify the proofs of  classical results. For simplicity of presentation, the method will be presented in the one--dimensional case.

\subsection{Paley--Wiener and Benrstein Spaces}

Given a bounded set $S\subset\R$, we denote by $PW_S$ the space of $L^2-$functions $f$ with spectrum in $S$. In other words, $PW_S$ consists of the functions $f$ such that
$$f(x)=\hat F(x):=\frac{1}{\sqrt{2\pi}}\int_S e^{-itx}F(t)\, dt, \ \ \ F\in L^2(S).
$$
Equipped with the $L^2-$norm $\|\cdot\|_2$, $PW_S$  is a Hilbert space.

 Given a compact set $S\subset\R$, we denote by $B_S$ the space of continuous bounded functions $f$ whose spectrum (in the sense of distributions) belongs to $S$. The latter means that
$$\int_\R f(x)\varphi(x)\, dx=0,$$
for every $\varphi$ from the Schwartz space whose inverse Fourier transform vanishes in some
neighborhood of $S$.
Equipped with the $L^\infty-$norm $\|\cdot\|_\infty$, $B_S$  is a Banach space.

\subsection{Sampling and Interpolation Sets for $PW_S$ and $B_S$}

 Let $\La$ be a u.d. set. It is called a set of stable sampling  (SS) for $PW_S$ if there is a constant $C$ such that
\begin{equation}\label{sam}
 \|f\|^2_2\leq C\sum_{\la\in\La}|f(\la)|^2, \ \ \forall f\in PW_S.
 \end{equation}
 $\La$ is called a set of interpolation (IS) for $PW_S$ if for every discrete function
$\{c_\la\}\in  l^2(\La)$ there exists $f\in PW_S$ such that
\begin{equation}\label{int}  f(\la) = c_\la, \ \ \la\in\La.\end{equation}

Similarly, $\La$ is called an SS for $B_S$, if if there is a constant $C$ such that
 $$
 \|f\|_\infty\leq C\sup_{\la\in\La}|f(\la)|, \ \ \forall f\in B_S.
 $$
 $\La$ is called an IS for $B_S$ if for every discrete function
$\{c_\la\}\in  l^\infty(\La)$ there exists $f\in B_S$ satisfying (\ref{int} ).

\section{Classical Results}

\subsection{Uniform Densities} The uniform densities of $\La$ play an important role for sampling and interpolation.
The lower  uniform density of a  u.d. set $\La$ is defined as
                   $$
 D^-(\Lambda):=\lim_{l\to\infty}\inf_{a\in{\R}}\frac{\#(\Lambda\cap(a,a+l))}{l},$$
 and its upper uniform density is
   $$
 D^+(\Lambda):=\lim_{l\to\infty}\sup_{a\in{\R}}\frac{\#(\Lambda\cap(a,a+l))}{l}.$$

\subsection{Sampling and Interpolation Theorems}
Assume $S=[a,b]$ is a single interval. Since the translations of the spectrum do not change the sampling/interpolation peroperty of $\La$, we may assume that $S=[-\sigma,\sigma]$. Denote by $PW_\sigma$ and $B_\sigma$ the Paley--Wiener and Bernstein spaces with the spectrum $S=[-\sigma,\sigma]$.

The Bernstein space $B_\sigma$ admits a complete characterization
of interpolation and sampling sets in terms of the uniform densities:

\begin{theorem} Let $\La$ be a u.d. set.

(i) ([Be89a])  $\La$ is an SS for $B_\sigma$  if and only if $D^-(\La)>\sigma/\pi.$

(ii) ([Be89b]) $\La$ is an IS for $B_\sigma$ if and only if $D^+(\La)<\sigma/\pi.$

\end{theorem}

The Paley--Wiener case is more involved: the interpolation  and sampling sets for $PW_\sigma$ can only be "essentially described" in terms
of the uniform densities:

\begin{theorem} Let $\La$ be a u.d. set.

(i) ([Be89a]) If  $D^-(\La)>\sigma/\pi$ then $\La$ is an SS for $PW_\sigma$.  If  $D^-(\La)<\sigma/\pi$ then $\La$ is not an SS for $PW_\sigma$.

(ii) ([Kh57])  If $D^+(\La) < \sigma/\pi$ then $\La$ is an IS for $PW_\sigma$.  If $D^+(\La) > \sigma/\pi$ then $\La$ is not an IS for $PW_\sigma$.

\end{theorem}

When the spectrum $S$ is a disconnected set, already when $S$ is a union of two disjoint intervals, simple examples show that no sharp sufficient condition for sampling and interpolation can be formulated in terms of some density of $\La$.
However, H. Landau discovered that the necessity of the density conditions in Theorems 1 and 2
is still valid in the general case:

\begin{theorem}
 ([La67]). Let $\La$ be a u.d. set and $S$ a bounded set.

(i) If $\La$ is an SS for $PW_S$ then $D^-(\La)\geq mes(S)/\pi$;

(ii) If $\La$ is an IS for $PW_S$ then $D^+(\La)\leq mes(S)/\pi$.

\end{theorem}

Here  $mes(S)$ denotes the measure of $S$.

Landau's theorem holds also  in several dimensions.

\section{The Approach}

The original proofs of the results above are quite involved. A  simpler proof of Theorem~3 is presented in [NO12], which also puts the results in a more general context.
 Beurling deduced part (i) of Theorem 2 from part (i) of Theorem 1 by using the so-called linear balayage operator ([Be89a], p. 349). A  simpler proof is presented in [OU12], see Theorem~9 below.


We prove

\medskip\noindent{\bf Claim}.  {\sl  The parts (i) and (ii) in Theorem 1 are equivalent in the sense
                     that each one can be deduced form the other one.                       The same is true for Theorems 2 and 3.}

\medskip

       This  claim saves "half of the job" in proving of Theorems 1--3.
                Its proof is based on an approach which is based on the three ingredients: the property of stability of sampling and interpolation, the complementarity property of the uniform densities and finally on the equivalence between sampling and interpolation in the special case of integer--sampling. The first two steps are similar for the Bernstein and Paley--Wiener cases. In the last step we use Theorem 8 for the Bernstein spaces, which is  more involved than Corollary 7 used for  the Paley--Wiener spaces.

\subsection{Frames and Riesz Sequences in Hilbert Spaces}We start with a duality result between frames and Riesz sequences.

Recall that a sequence of vectors $u_n$ in $H$ is a  frame if there exist numbers $A,B >0$ such that for all $x\in H$ we have
\begin{equation}\label{frame}A\Vert x\Vert_H^2\leq\sum_{n}|\langle x,u_n\rangle|^2\leq B\Vert x\Vert_H^2.\end{equation}

 A family of vectors $\{u_n\}$ is called                              a  Riesz sequence in $H$ if there are
                            constants $ c,C>0$ such that
\begin{equation}\label{rs}
                C \Vert\{ c_n\}\Vert_{l^2}   \geq       \Vert\sum_n c_nu_n\Vert_H\geq c \Vert\{ c_n\}\Vert_{l^2},
\end{equation} for every finite sequence $\{c_n\}$.


\begin{proposition}\label{p3}
Assume that a set $U$ is an orthonormal basis in a Hilbert space $H$. Assume
$H$ is a direct sum of two orthogonal subspaces $H_1$ and $H_2$, and denote by $P_j$ the
orthogonal projections on $H_j$. Assume further that $U$ is a union of two disjoint
subsets $V$ and $W$. Then the following statements are equivalent:

(i) $P_1V$ is a frame in $H_1$;

(ii) $P_2W$ is an Riesz sequence in $H_2$.
\end{proposition}

        \medskip\noindent  {\bf   Proof}. (i) $\Rightarrow$ (ii).
Take any $l^2-$sequence $c(w)$ and write
$$
f=\sum_{w\in W}c(w)w=:f_1+f_2, \ \ f_1\in H_1, f_2\in H_2.
$$
For every $v\in V$ we have  $\langle f, v\rangle=0$, and so $\langle f_1, v\rangle=-\langle f_2,v\rangle$. This and property (i) imply
$$
\| f_1\|_H^2\leq\frac{1}{A}\sum_{v\in V}|\langle f_1,v\rangle|^2=\frac{1}{A}\sum_{v\in V}|\langle f_2,v\rangle|^2$$$$\leq\frac{1}{A}\sum_{u\in U}|\langle f_2,u\rangle|^2=\frac{1}{A}\| f_2\|_H^2,
$$
where $A$ is the constant from (\ref{frame}).
We conclude that
$$
\sum_{w\in W} |c(w)|^2=\Vert f\Vert_H^2=\|f_1\|_H^2+\|f_2\|_H^2\leq (1+\frac{1}{A})\|f_2\|_H^2,
$$which proves the right inequality in (\ref{rs}). The left inequality is obvious.

The proof of (ii) $\Rightarrow$ (i) is similar.

\subsection{Stability of Sampling and Interpolation}
By stability we mean that a slight perturbation of a sampling (interpolation) set remains to be a sampling (interpolation) set.

Given $\delta>0$ and a u.d. set $\La$, we call $\La'$ a $\delta-$perturbation of $\La$ if $\La'$ is a u.d. set which  admits a representation $\La'=\{\la+\epsilon_\la, \la\in \La\}$, where $ \sup_{\la\in\La}|\epsilon_\la|\leq\delta.$
\begin{proposition}\label{p1} Suppose  $\La$ is a u.d. set.

(i) If $\La$ is an SS for $PW_S$, then there exists $\delta>0$ such that every  $\delta-$perturbation of $\La$ is an SS for $PW_S$.

(ii) A similar property holds for the interpolation sets.\end{proposition}

A similar to Proposition \ref{p1} result is true for the Bernstein spaces $B_S$.

These results are well--known, see [Yo01, ch.4, Theorems 11 and 13]. For the readers convenience, we show how they can be deduced with the help of the "restriction operator" $$ R:f\to f|_\La. $$


Observe that the inequality converse to (2) holds (see [Yo01, ch.4, Theorem 4]): Given a u.d. set $\La$, there exists $C$ such that \begin{equation}\label{bessel} \sum_{\la\in\La}|f(\la)|^2\leq C\|f\|_2, \ \ \ \forall f\in PW_S.\end{equation}The constant $C$ above depends only on the diameter of $S$ and the magnitude of  separation constant $d(\La)$ defined in (1).

Inequality (\ref{bessel})  shows that the restriction operator $R$ is a bounded linear operator from $PW_S$ to $l^2(\La)$.
Clearly, we have:

\medskip

(a) $\La$ is an IS for $PW_S$ if and only if the restriction operator is surjective;

(b) $\La$ is an SS for $PW_S$ if and only if there exists $C$ such that $$\|f\|_2\leq C\|Rf\|_{l^2}, \ \ \ \forall f\in PW_S.$$

Let $\La'=\{\lambda+\epsilon_\la,\la\in\La\}$ where $|\epsilon_\la|\leq \delta, \la\in\La$, be a $\delta-$perturbation of $\La$, and let
$R'f=f|_{\La'}$ be the corresponding restriction operator. We may consider  $R'$ as an operator from $PW_S$ to $l^2(\La)$. Observe that
 $$
|f(\la+\epsilon_\la)-f(\la)|=|\int_\la^{\la+\epsilon_\la}f'(x)\,dx|\leq |f'(\xi_\la)|\epsilon_\la\leq\delta|f'(\xi_\la)|,
$$with some $\xi_\la$ between $\la$ and $\la+\epsilon_\la$. 

 Assume that $\delta<d(\La)/4$ and that $S\subset[-\si,\si],$ for some $\si>0$. Denote by $F$ the inverse Fourier transform of $f$. Clearly, the sequence $\{\xi_\la,\la\in\La\}$ can be represented as the union of two u.d. sequences whose separation constants are greater than $d(\La)$.  We now apply (\ref{bessel}) to $f'$  to obtain
$$
\Vert (R'-R)f\Vert_{l^2}^2=\sum_{\la}|f(\la+\epsilon_\la)-f(\la)|^2\leq \delta^2\sum_{\la}|f'(\xi_\la)|^2$$ $$\leq C\delta^2\Vert f'\Vert^2_2=C\delta^2\int_{-\si}^{\si}
|tF(t)|^2\,dt\leq C\si^2\delta^2\Vert f\Vert^2_{2}.
$$
We see that $\Vert R-R'\Vert$ can be made arbitrarily small, provided $\delta$ is sufficiently small.

The stability of  sampling property is now evident: One just needs to choose $\delta$ so small that $\Vert R-R'\Vert$ is less than the constant $C$ is the property (b) above.

To prove the stability of interpolation, one may use (a) and the well-known fact on perturbations of a surjective linear operator: {\sl Assume $X,Y$ are Banach spaces, and a bounded linear operator $A:X\to Y$ is a surjection. Then there exists $\gamma=\gamma(A)>0$ such that for every linear operator $B:X\to Y$ with $\Vert B\Vert<\gamma,$ the operator $A+B:X\to Y$ is a surjection.}

The proof of stability of sampling and interpolation for the Bernstein spaces is  similar.

\subsection{Complementarity of $D^-$ and $D^+$}
The upper and lower uniform densities are "complementary" in the following sense:
\begin{proposition}\label{p2}Suppose $\La\subset\delta\Z$. Then
$$D^-(\La)+D^+(\delta\Z\setminus\La)=\frac{1}{\delta}.$$
\end{proposition}The proof is trivial, and we omit it.

\subsection{Integer Sampling. Paley--Wiener Spaces} Here we apply Proposition 4 to the case $H=L^2(\T), H_1=L^2(S), H_2= L^2(\T\setminus S)$ and $U=\{e^{i n t}, n\in\Z\}$,
where $\T=[0,2\pi]$ is the unite circle with the  normalized Lebesgue measure $dt/2\pi$:

\begin{corollary}\label{c1}
Suppose $S\subset \T$ and $\La\subset\Z$. Then the following statements are equivalent:

(i) $\La$ is an SS for $PW_S$;

(ii)  $\Z\setminus \La$ is an IS for $PW_{\T\setminus S}$.
\end{corollary}

To deduce this corollary from Proposition 5 it suffices to use the following well-known  correspondence (see [Yo01, ch.4]):

\medskip

 (c) A u.d. set $\La$ is an SS for $PW_S$ if and only if the exponential system $\{e^{i\la t},\la\in\La\}$ is a frame in $L^2(S)$;

  (d) A u.d. set $\La$ is an IS for $PW_S$ if and only if the exponential system $\{e^{i\la t},\la\in\La\}$ is a Riesz sequence in $L^2(S)$.

  \medskip

  More generally, suppose that  $S\subset(0, 2\pi/\delta)$ and $\La\subset\delta\Z$. Using a linear change of variables, one may deduce from Corollary 7 that  $\La$ is an SS for $PW_S$ if and only if $\delta\Z\setminus \La$ is an IS for $PW_{(0,2\pi/\delta)\setminus S}$.

\subsection{Integer Sampling. Bernstein Spaces}

\begin{theorem}\label{t5}
 Let $\La\subset\Z,  0< a<2\pi$. Then the following statements are equivalent:

 (i) $\La$ is an IS for $B_{[0,a]}$;

 (ii) $\Z\setminus\La$ is an SS  for $B_{[a,2\pi]}$.
\end{theorem}

This result  is more difficult than Corollary \ref{c1}. We do not know if an analogue of Theorem~\ref{t5} remains true for the Bernstein space $B_{S}$  even when $S$ is a finite union of  intervals on $[0, 2\pi]$.

More generally, suppose that $\delta>0,$ $0< a<2\pi/\delta$ and  $\La\subset\delta\Z$. It follows from Theorem~\ref{t5} that $\La$ is an SS for $B_{[0,a]}$ if and only if $\delta\Z\setminus\La$ is an SS  for $B_{[a,2\pi/\delta]}$.

We will prove Theorem \ref{t5} in sec. 5.

\section{Proof of Claim} Here we prove the Claim.   In the case of Theorem 1 we use Theorem 8, while in the case of Theorems 2 and 3 we use Corollary 7. The rest of the proof is the same. So, for simplicity, below we present the method in the case of  Theorem 3.

Let us assume that  condition (i) of Theorem 3 holds true and prove  condition (ii).

Fix a bounded set  $S\subset\R$ and assume that a u.d. set $\La$ is an IS for $PW_S$. We may assume that $S\subset(0,\infty)$.

We consider two cases:

(e) Let us assume additionally that $\La$ is a subset of $\delta\Z$, for some $\delta> 0$. It follows from Corollary \ref{c1} that the set $\Gamma:=\delta\Z\setminus \La$ is an SS for $PW_G,$ where $G=(0, 2\pi/\delta)\setminus S$. By (i), $D^-(\Gamma)\geq mes(G)/2\pi=1/\delta-mes(S)/2\pi$. Using Proposition 6, we conclude that
$$
D^+(\La)=\frac{1}{\delta}-D^-(\Gamma)\leq\frac{mes(S)}{2\pi}.
$$

¤

(f) Assume $\La$ is not a subset of an arithmetic progression. Choose any set $\La'\subset\delta\Z$ which is a $\delta-$perturbation of $\La$.  From the definition of the upper uniform density, it is evident that $D^+(\La')=D^+(\La)$. By Proposition 5, we may assume that $\delta$ is so small that $\La'$ is an IS for $PW_S$. Hence, by the previous step, we have
$D^+(\La)=D^+(\La')\leq mes(S)/2\pi$.

The proof of (ii) $\Rightarrow$ (i) is similar.

\section{Proof of Theorem \ref{t5}}

\subsection{Connection between Sampling/Interpolation in Paley--Wiener and Bernstein Spaces}
The following theorem is proved in [OU12, Theorem 2.1]:

\begin{theorem}  Given a compact $S \subset \R$, a u.d. set $\La$ and $\delta > 0$.

(i) If $\La$ is an SS for $B_{S+[-\delta,\delta]}$, then $\La$ is an SS for $PW_S$.

(ii) If $\La$ is an SS for $PW_{S+[-\delta,\delta]}$, then $\La$ is an SS for $B_S$.

(iii) If $\La$ is an IS for $B_{S}$, then $\La$ is an IS for $PW_{S+[-\delta,\delta]}$.

(iv) If $\La$ is an IS for $PW_{S}$, then $\La$ is an IS for $B_{S+[-\delta,\delta]}$.

\end{theorem}


The proof of this theorem in [OU12] is fairly elementary.

A similar result holds in several dimensions, see [OU12].

\subsection{Weak Limits of Translates}
Let  $\La_j$ be u.d. sets in $\R$ satisfying $$d(\La_j) > d >0,$$ where $d(\La_j)$ are the separation constants defined in (\ref{sepc}).
   A set   $\La\subset\R$ is called the weak limit of $\La_j$ as $j\to\infty$, if for every $\epsilon>0$
                       and for every interval $I=(a,b), a,b\not\in\La$,  the set $\La\cap I$ is an $\epsilon-$perturbation of $\La_j\cap I$ for all but a finite number of $j$'s.

  Denote by $W(\La)$ the set of all weak limits of the translates of $\La$: a set $\La'$ belongs to $W(\La)$ if there is a sequence of translates $\La-a_j$  which converges weakly to $\La'$.

  Weak limits preserve the sampling and interpolation properties:

  \begin{proposition}\label{pp} Suppose $\La$ is a u.d. set.

  (i)  ([Be89a], p. 343) If $\La$ is an SS for $B_\si$ then every $\La'\in W(\La)$ is an SS for $B_\si$.

  (i)  ([Be89b], p. 351) If $\La$ is an IS for $B_\si$ then every $\La'\in W(\La)$ is an IS for $B_\si$.

  \end{proposition}

Let us say that a set $\La$ is a uniqueness set (US) for $B_\si$, if no non-trivial function from $B_\si$ vanishes on $\La$.
There is a nice connection between sampling and uniqueness sets discovered by Beurling:

\begin{theorem}\label{tb} ([Be89a, p. 345]) A u.d. set $\La$ is an SS for $B_\si$ if and only if  every $\La' \in W(\La)$ is a US for $B_\si$.
\end{theorem}

\begin{corollary}\label{c2}  If $\La$ is an SS for $B_\si$ then it is also an SS for $B_{\si+\epsilon}$, for all sufficiently
small $\epsilon>0$.
\end{corollary}

\noindent{\bf Proof} (by contradiction). Assume that this statement is not true. By Theorem \ref{tb}, for every $j \in\N$ there is a function $f_j\in B_{\si+1/j}$ with $\|f_j\|_\infty=1$ which vanishes on some $\La_j\in W(\La)$. We may assume that $|f(0)|>1-1/j$ (otherwise, we consider an appropriate translation of $\La_j$).  We may choose a subsequence $j_k$ such that $\La_{j_k}$ converge weakly to some set $\Gamma$. One may easily check that $\Gamma\in W(\La)$. Further, it is clear that  $\{f_j\}$ is a normal family. Choosing a subsequence of $j_k$, we get a non-trivial limiting
function $f\in B_\si$ which vanishes on  $ \Gamma$. By Theorem \ref{tb}, this means that $\La$ is not an SS for $B_\si$.

\subsection{Auxiliary Lemmas}

\begin{lemma}\label{l1}   Suppose $\La$ is an IS for $B_\sigma.$ Then $\La$ is not a US for $B_\sigma$.\end{lemma}

\noindent{\bf Proof} (by contradiction). If  $\La$ is both an IS and a US for $B_\sigma$, then the restriction operator
      $R: B_\sigma\to l^\infty$ is a  bounded bijection. Then the inverse operator also is bounded, so that $\Vert f\Vert_\infty\leq C\Vert Rf\Vert_{l^\infty}$. This means that $\La$ is an SS for $B_\sigma$.

  Fix any $\la_0\in\La$  and let   $f\in B_\sigma$ such that $f(\lambda_0)=1$ and $f$  vanishes on $\La\setminus\{\lambda_0\}$. The function $f(x)\sin\epsilon(x-\la_0)$ belongs to $B_{\si+\epsilon}$ and vanishes on $\La$, so that $\La$ is not an SS for $B_{\si+\epsilon}$. This contradicts to  Corollary \ref{c2}.

\begin{lemma}\label{l2}  Suppose $\La$ is an IS for $B_\sigma.$ Then for every point  $l\in\R\setminus\La$ the set $\La\cup\{l\}$ is an IS for $B_\sigma$.\end{lemma}

\noindent{\bf Proof}. Fix any point $l\in\R\setminus\La$. By Lemma \ref{l1} , there exists $\varphi\in B_\sigma$ which vanishes on $\La$. We may assume  that $\varphi(l)\ne0$ (otherwise we consider the function $\varphi(x)/(x-l)^n$, where $n$ is the order of zero of $\varphi$ at $l$). We may also assume that $\Vert \varphi\Vert_\infty=1$.

    Given any data $c\in l^\infty(\La\cup\{l\})$, let $f\in B_\sigma$ be such that $f(\la)=c_\la,\la\in\La.$
      Set
     $$
     g(x):=f(x)+\frac{c_l-f(l)}{\varphi(l)}\varphi(x).
     $$
     Clearly, $g(\lambda)=c_\lambda,\lambda\in\La$ and $g(l)=c_l.$ This shows that $\La\cup\{l\}$ is an IS for $B_\sigma$.

\begin{lemma}\label{l3} Suppose $f\in B_{[0,2\pi]}$ and $|f(x)|=O(x^{-2}), |x|\to\infty$. Then $$
\sum_{n\in\Z}f(n)=0.
$$
\end{lemma} To prove this, one may consider the Fourier transform $F=\hat f$. It is easy to see that its Fourier coefficients  belong to $l^1$. Then one may just write its Fourier  series at the origin.

\begin{lemma}\label{l4} Suppose $f\in B_\si, 0<\si<\pi$. If $f$ vanishes on $\Z$ then $f=0$.
\end{lemma}

This follows from the Jensen's Formula ([Le96], Lec.2].

      \subsection{Proof of Theorem \ref{t5}}
       1.  Suppose $\La$ is an IS for $B_{[0, a]}$. We must show that the set $Q:=\Z\setminus\La$ is an SS for $B_{[a,2\pi]}$. Suppose it is not. Then, by Theorem \ref{tb}, there is a weak limit $Q'$ of translates of $Q$,  which is not a uniqueness set for $B_{[a,2\pi]}$.  Taking if necessary an additional translation, we may assume that $Q'\subset\Z$.

       It is clear that $\La'=\Z\setminus Q'$ is a weak limit of translates of $\La$. So, by Proposition \ref{pp},  $\La'$ is an IS for $B_{[0,a]}$.

        Thus, we get a partition $\Z=\La'\cup Q'$ such that

\medskip
        (g)  $\La'$ is an IS for $B_{[0,a]}$;

        (h) $Q'$ is not a uniqueness set for $B_{[a,2\pi]}$.

        \medskip

        We show that this leads to a contradiction.

        By (h), we may fix a non-trivial function $h\in B_{[a,2\pi]}$ which vanishes on $Q'$. This function cannot vanish on $\Z$, due to Lemma \ref{l4}. Fix a point $x_0\in\La'$ such that $h(\lambda_0)\ne0$.

         Take any point $x_1\not\in\La'$. Now we use  (i). By Lemmas \ref{l1} and \ref{l2}, there is a non-trivial function $\varphi\in B_{[0,a]}$ satisfying $\varphi|_{\La'\cup\{x_1\}}=0.$

Set
       $$
       g(x):=\frac{\varphi(x)}{(x-x_1)(x-x_0)^n}\in B_{[0,a]},
       $$where $n$ is the multiplicity of the zero of $\varphi$ at $x_0$. Then $g(x_0)\ne0$, $g|_{\La'\setminus\{x_0\}}=0$ and
       $|g_\lambda(x)|=O(x^{-2}),|x|\to\infty.$

Finally, consider the function $f=gh$. Clearly, $f$ satisfies he assumptions of  Lemma \ref{l4}. Hence,
  $$ 0= \sum_{n\in\Z}f(n)=   g(x_0)h(x_0)\ne0.
       $$This is a contradiction.


\medskip\noindent 2.        Suppose $\Z\setminus\La$ is an SS  for $B_{[a,2\pi]}.$ By Corollary \ref{c2},  there exists $\epsilon>0$ such that
      $\Z\setminus\La$ is an SS  for $B_{[a-2\epsilon,2\pi]}.$ Then, by Theorem 6 (i), $\Z\setminus\La$ is an SS  for $PW_{[a-\epsilon,2\pi]}.$ By Corollay 7, $\La$ is an IS for $PW_{[0, a-\epsilon]}$. By Theorem 6 (iv), we conclude that $\La$ is an IS for $B_{[0, a]}$.


\begin{thebibliography}{99}




\bibitem[Be89a]{B} Beurling, A. Balayage of Fourier--Stiltjes Transforms. In:
The collected Works of Arne Beurling, vol.2, Harmonic Analysis.
 Birkhauser, Boston, 1989.


\bibitem[Be89b]{Bi} Beurling, A. Interpolation for an interval in $\R$. In: The Collected Works of Arne Beurling, in: Harmonic Analysis, vol. 2, Harmonic Analysis. Birkh\"{a}user, Boston,
1989.



\bibitem[Kh57]{K} Kahane, J.-P. Sur les fonctions
moyenne-p\'{e}riodiques born\'{e}es.  Ann. Inst. Fourier, 7, (1957),
293--314.

 \bibitem[La67]{L2} Landau, H. J. Necessary density conditions for
   sampling and interpolation of certain entire functions. Acta Math. 117, (1967), 37--52.

   \bibitem[Le96]{Levin} Levin, B. Ya.  Lectures on entire functions. In collaboration with and with a preface by Yu. Lyubarskii, M. Sodin and V. Tkachenko. Translated from the Russian manuscript by Tkachenko. Translations of Mathematical Monographs, 150. American Mathematical Society, Providence, RI, 1996.

\bibitem[NO12]{NO2}   Nitzan, S., Olevskii, A. Revisiting Landau's density theorems for Paley--Wiener spaces. C. R. Mathematique  350, no. 9--10, (2012), 509--512.



   \bibitem[OU12]{ou12} Olevskii, A., Ulanovskii, A. On multi-dimensional sampling and interpolation. Anal. Math. Phys. 2 (2012), no. 2, $149--170.$


\bibitem[Yo01]{Ya} Young, R.M. An introduction to Nonharmonic Fourier Series. Academic Press. 2001.
\end{thebibliography}
\end{document}